# Machine Learning based parameter tuning strategy for MMC based topology optimization

**Xinchao Jiang[a,*], Hu Wang[a,**], Yu Li[a], Kangjia Mo[a]**

*a. State Key Laboratory of Advanced Design and Manufacturing for Vehicle Body, Hunan University, Changsha, 410082, P.R. China*

## Abstract

Moving Morphable Component (MMC) based topology optimization approach is an explicit algorithm since the boundary of the entity explicitly described by its functions. Compared with other pixel or node point-based algorithms, it is optimized through the parameter optimization of a Topological Description Function (TDF). However, the optimized results partly depend on the selection of related parameters of Method of Moving Asymptote (MMA), which is the optimizer of MMC based topology optimization. Practically, these parameters are tuned according to the experience and the feasible solution might not be easily obtained, even the solution might be infeasible due to improper parameter setting. In order to address these issues, a Machine Learning (ML) based parameter tuning strategy is proposed in this study. An Extra-Trees (ET) based image classifier is integrated to the optimization framework, and combined with Particle Swarm Optimization (PSO) algorithm to form a closed loop. It makes the optimization process be free from the manual parameter adjustment and the reasonable solution in the design domain is obtained. In this study, two classical cases are presented to demonstrate the efficiency of the proposed approach.

---

[*] First author. *E-mail address*: jiangxinchao@hnu.edu.cn (X.C. Jiang)

[**] Corresponding author. Tel.: +86 0731 88655012; fax: +86 0731 88822051.

*E-mail address:* wanghu@hnu.edu.cn (H. Wang)

## 1. Introduction

Topology optimization [1, 2], a branch of design optimization, is a mathematical method to solve a material layout problem constrained to a given design domain, loading, and boundary conditions. This method determines the optimal distribution of material and the corresponding structure which satisfies the desired properties (e.g. compliance) under the design constraints. Since the pioneering work of Bendsoe and Kikuchi [3], topology optimization has received considerable research attention. There are numerous topology methods such as Solid Isotropic Material Penalty (SIMP) [4-6], Level Set Method (LSM) [7, 8], evolutionary method [9] and MMC based method [10-12]. These methods have been applied successfully to a wide range of topological design problems. For example, SIMP method has been integrated into commercial software such as Abaqus [13] due to its advantages of simplified program, fast calculation and easy implementation. LSM can be used for both shape optimization and topology optimization, its results are smooth boundary with no checkerboard and gray scale elements [14], etc. However, most approaches are considered to do topology optimization in an implicit way, while the MMC based method is a more explicit and geometrical way. As for SIMP, it means the optimal structural topology is identified from a black-and-white pixel image. For LSM, it means the optimum can be identified from the level set of a TDF defined in a prescribed design domain. There are some problems associated with the implicit approached. Firstly, under the framework of implicit topology optimization, it is difficult to give the precise control of structural feature size, which is very important from the perspective of manufacture [15]. It results from the lack of explicit geometry information in the implicit optimization model [16]. Hence, implicit optimization models are difficult to establish a direct link between optimization models and Computer-Aided-Design (CAD) modeling systems. Secondly, the number of design variables involved in implicit topology optimization approaches is connected with the resolution since they are developed within the pixel (SIMP) or node point based (LSM) solution framework. It is relatively large especially for 3-dimentional problems. etc.

In order to address these issues, a so-called MMC based topology optimization framework is established in Refs [10-12]. This approach considers a set of components as the optimization unit and achieves the structural topology optimization by morphing, merging, and overlapping operations between the components. It is controlled by the parameters of TDF of each component and it does not concern the finite element grid division. Fig. 1 shows the basic idea of this approach schematically. Since this approach doing topology optimization in a more explicit and geometrical way, it has flourished in extensive applications, e.g. three-dimensional problems [17, 18], geometric size control [19], additive manufacturing, and stress constraints [20], etc.

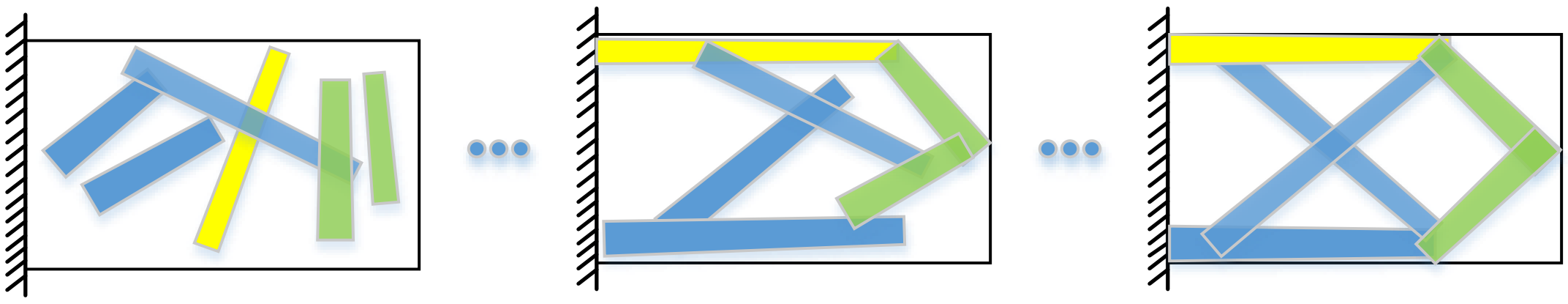

**Fig. 1** The basic idea of the MMC-based topology optimization approach

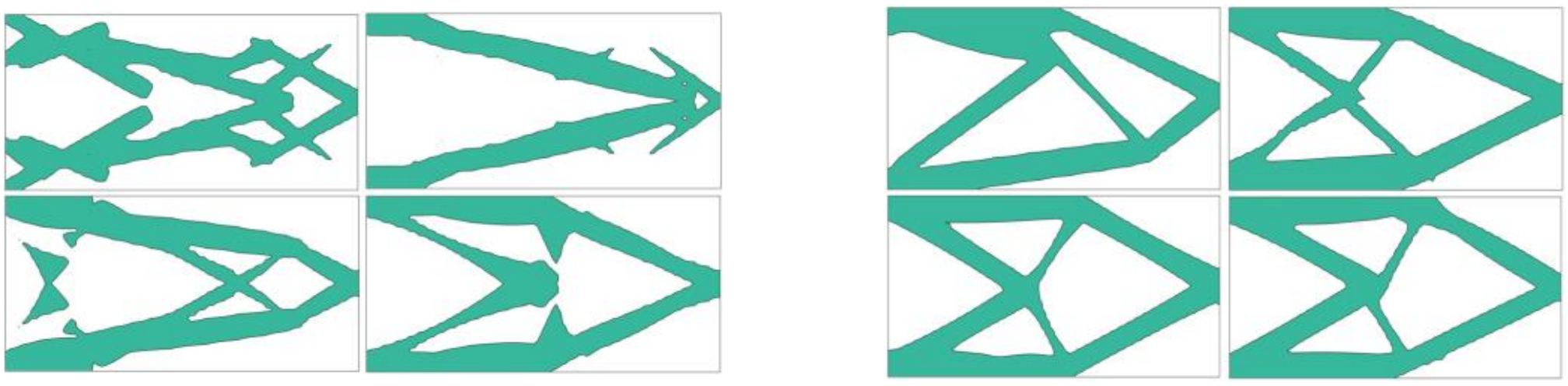

**(a)** **(b)**

**Fig. 2** Short beam optimization results (a) infeasible results (b) feasible result

Although it has many distinctive advantages, its performance mainly relies on experienced parameters. In MMC based method, MMA is employed as the optimizer [21]. It is might be time-consuming to tune the parameters of MMA. Furthermore, as shown in the Fig. 2, it might obtain infeasible solutions with inappropriate parameter setting. At the same time, the feasible solution might not be easily obtained according to experience. With the attempt to address these bottlenecks and make it be free from the manual parameter adjustment. This study proposed a ML based parameter tuning strategy for the MMC based approach. Actually, due to the state-of-art performance of ML, it has been extended in the field of topology optimization already. Sharad [22] used Convolutional Neural Network (CNN)-based predictor to accelerate the SIMP-based topology optimization. Compared with Sharad who used supervised learning to predict topology optimization, Sosnovik [23] employed unsupervised learning. He combined CNN and Generative Adversarial Networks (GANs) to predict the optimized structure, which is a complementary alternative to the traditional topology optimization. Banga [24] proposed a deep learning approach based on a 3D encoder-decoder Convolutional Neural Network architecture for accelerating 3D topology optimization and to determine the optimal computational strategy for its deployment. Moreover, Supported Vector Regression (SVR) as well as K-Nearest-Neighbors (KNN) ML models were employed to the MMC by Lei [25] to obtain the optimized distribution in a prescribed design domain almost instantaneously once the objective/constraint functions and external stimuli/boundary conditions are specified. Compared with these existing ML based topology optimization methods. This study focuses on parameter tuning strategy to guarantee parameter selection depends on experience less in MMC based approach. Two numerical cases are presented to demonstrate the efficiency of this strategy.

The remainder of the paper is organized as follows. In Section 2, the structural topology optimization under the MMC based framework and corresponding formulations are introduced. Section 3 is devoted to establishing the framework of parameter tuning strategy and the involved tools are also discussed. In Section 4, the numerical cases are presented to demonstrate the efficiency of the proposed method. Finally, some conclusions are provided in Section 5.

## 2. MMC based topology optimization framework

### 2.1. Structural shape and topology description

In the MMC based approach [10], the structural topology descriptions about components can be constructed in the following way.

$$\begin{cases} \phi^s(\boldsymbol{x}) > 0, & \text{if } \boldsymbol{x} \in \Omega^s \\ \phi^s(\boldsymbol{x}) = 0, & \text{if } \boldsymbol{x} \in \partial\Omega^s \\ \phi^s(\boldsymbol{x}) < 0, & \text{if } \boldsymbol{x} \in \mathbb{D} \setminus \Omega^s \end{cases} \tag{1}$$

where $\mathbb{D}$ means a predefined design domain and $\Omega^s \subset \mathbb{D}$ represents a set of components $\phi^s(\boldsymbol{x}) = \max(\phi_1,...,\phi_n)$ [26]. These components represent the parts (i.e. entities) that have materials during the optimization process. Actually

$$\begin{cases} \phi_i(\boldsymbol{x}) > 0, & \text{if } \boldsymbol{x} \in \Omega_i \\ \phi_i(\boldsymbol{x}) = 0, & \text{if } \boldsymbol{x} \in \partial\Omega_i \\ \phi_i(\boldsymbol{x}) < 0, & \text{if } \boldsymbol{x} \in \mathbb{D} \setminus \Omega_i \end{cases} \tag{2}$$

where $\Omega_i$ is the region occupied by the $i$-th component and $\Omega^s = \bigcup_{i=1}^{n} \Omega_i$. Its geometry representation can be illustrated schematically in Fig. 3.

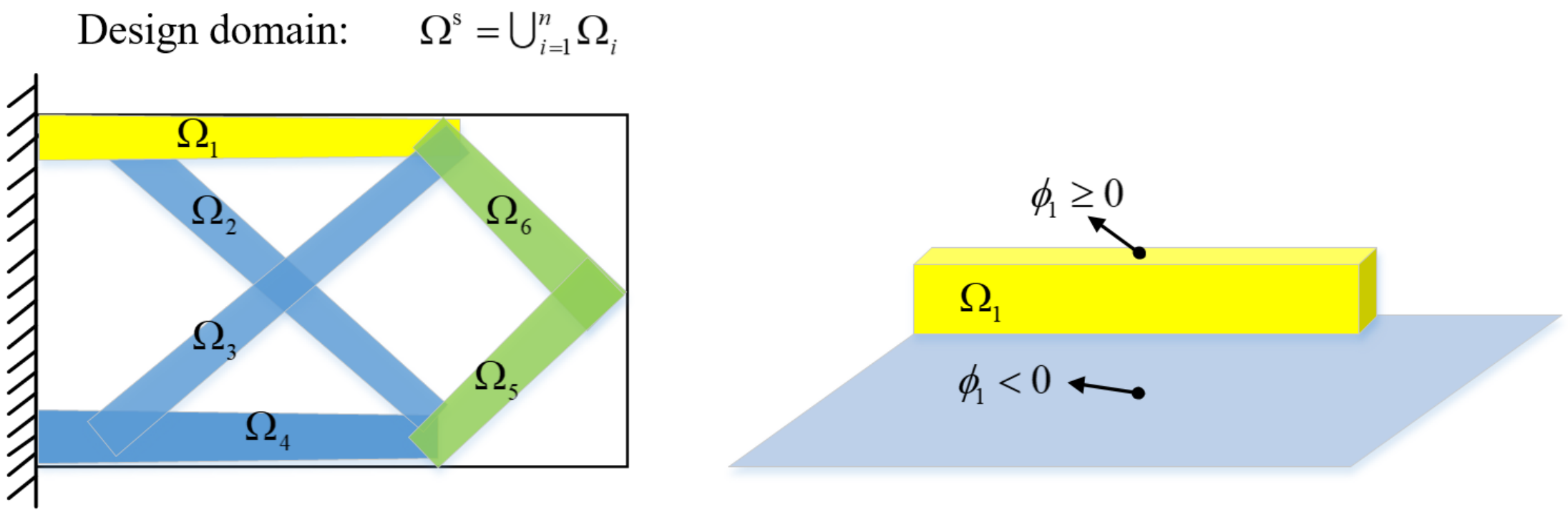

**Fig. 3** The representation of each components through the level set functions

The descriptive equation for each component is explicit [27, 28], it can be expressed directly in terms of certain function with $\boldsymbol{x}$ and $\boldsymbol{y}$.

$$\phi_i(\boldsymbol{x}, \boldsymbol{y}) = \left(\frac{\boldsymbol{x}'}{L_i}\right)^p + \left(\frac{\boldsymbol{y}'}{f(\boldsymbol{x}')}\right) - 1 \tag{3}$$

with

$$\begin{Bmatrix} \boldsymbol{x}' \\ \boldsymbol{y}' \end{Bmatrix} = \begin{bmatrix} \cos\theta_i & \sin\theta_i \\ -\sin\theta_i & \cos\theta_i \end{bmatrix} \begin{Bmatrix} \boldsymbol{x} - \boldsymbol{x}_{0i} \\ \boldsymbol{y} - \boldsymbol{y}_{0i} \end{Bmatrix} \tag{4}$$

$p$ is a relatively large even integer number (in this study, $p = 6$), $\boldsymbol{x}_{0i}, \boldsymbol{y}_{0i}$ denotes the coordinate of the center of the component, $L_i$ denotes the half length of the component and $\theta_i$ is the inclined angle of the component respectively. The quadratic varying thickness component with straight skeleton in the form of hyperelliptic equation is adopted in this study, which is presented in Fig. 4. The shape of a component is controlled by $f(\boldsymbol{x}')$. To sum up, the layout (i.e., shaper and topology) of a structure can be solely determined by a design vector $\boldsymbol{D}=((\boldsymbol{D}^1)^T,...,(\boldsymbol{D}^n)^T)^T$. Here $\boldsymbol{D}^i=(x_0,y_0,L_i,\theta_i,\boldsymbol{d}_i^T)^T$ and the symbol $\boldsymbol{d}_i$ denotes the vector of parameters associated with $f(\boldsymbol{x}')$ (e.g., $t_1$, $t_2$ and $t_3$ in Fig. 4).

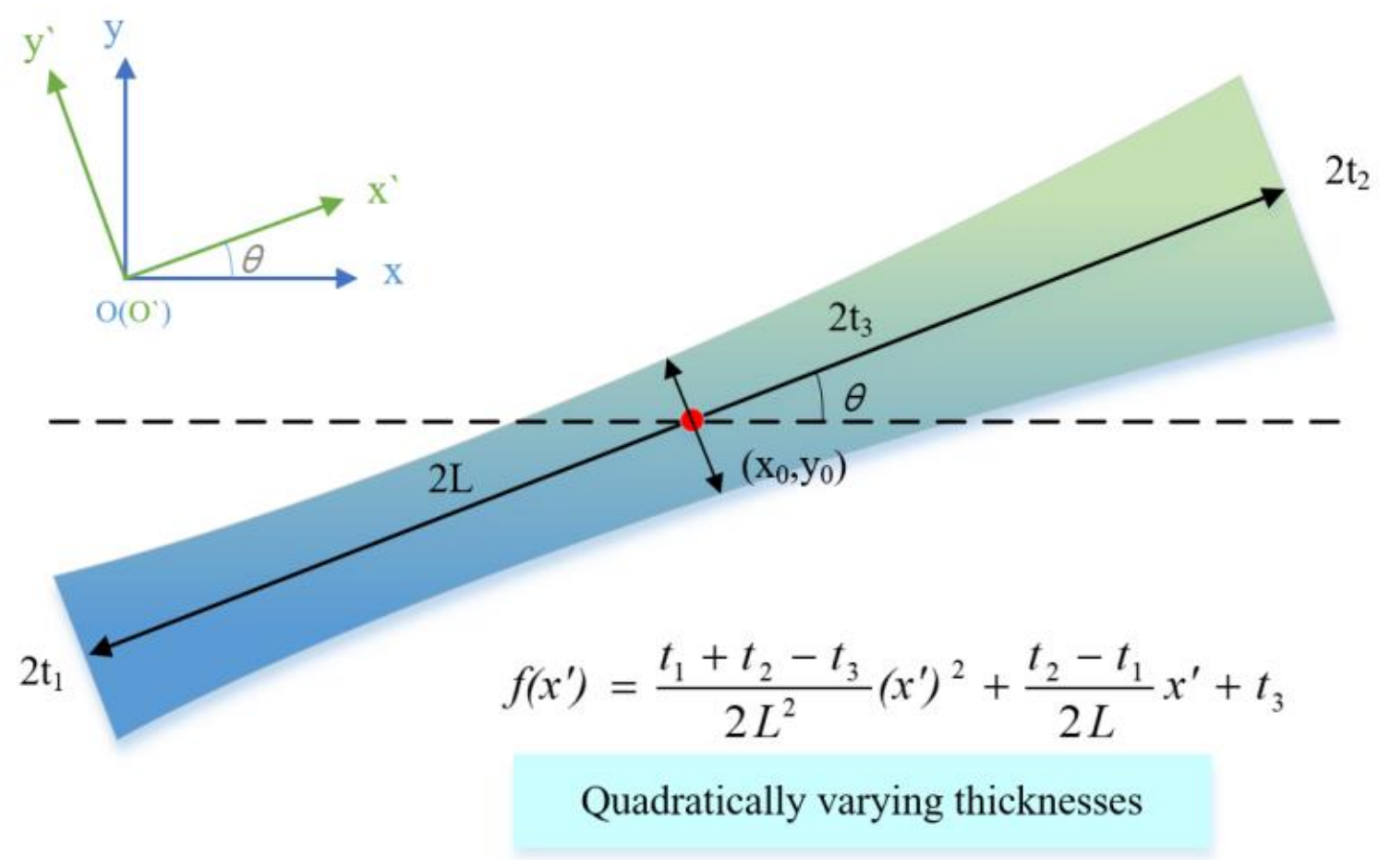

**Fig. 4** Geometry description of the variant thickness component

## 2.2. Problem formulation

In this study, compliance minimization with available volume constraint is considered. The mathematic model of the topology optimization under MMC based framework can be formulated as follows:

find

$$\boldsymbol{D}=\left((\boldsymbol{D}^1)^T,\ldots,(\boldsymbol{D}^i)^T,\ldots,(\boldsymbol{D}^n)^T\right)^T,\boldsymbol{u}(\boldsymbol{x}) \tag{5}$$

minimize

$$C=\int_D H\left(\phi^s\left(\boldsymbol{x};\boldsymbol{D}\right)\right)\boldsymbol{f}\bullet\boldsymbol{u}dV+\int_{\Gamma_t}\boldsymbol{t}\bullet\boldsymbol{u}dS \tag{6}$$

subject to

$$\int_D\left(H\left(\phi^s\left(\boldsymbol{x};\boldsymbol{D}\right)\right)\right)^q\mathbb{E}:\varepsilon(\boldsymbol{u}):\varepsilon(v)dV=\int_D H\left(\phi^s\left(\boldsymbol{x};\boldsymbol{D}\right)\right)\boldsymbol{f}\bullet vdV \\ +\int_{\Gamma_t}\boldsymbol{t}\bullet vdS,\quad\forall v\in u_{ad} \tag{7}$$

$$\int_D H\left(\phi^s\left(\boldsymbol{x};\boldsymbol{D}\right)\right)dV<\bar{V} \tag{8}$$

$$\boldsymbol{D} \subset u_D \tag{9}$$

$$\boldsymbol{u} = \bar{\boldsymbol{u}}, \text{on } \Gamma_u \tag{10}$$

$$u_{ad} = \{ v \mid v \in H^1(\boldsymbol{D}), v = 0 \text{ on } \Gamma_u \} \tag{11}$$

with

$$H = H(x) = \begin{cases} 0, & \text{if } x \le 0 \\ 1, & \text{otherwise} \end{cases} \tag{12}$$

From Eqs. (5) to (11), only a single-phase material problem is considered and the Galerkin numerical solution method is adopted, where $H(x)$ is Heaviside function [29], $\boldsymbol{D}^i$ is vector of design variables associated with the *i*-th component respectively. $H\left(\phi^s\left(\boldsymbol{x};\boldsymbol{D}\right)\right)$ converts $\phi^s\left(\boldsymbol{x};\boldsymbol{D}\right)$ into solid-void material distribution. Here $q$ is an integer larger than 1, in this study, $q=2$. $\boldsymbol{f}$ and $\boldsymbol{t}$ denote the body force density in $\Omega_i$ and the surface traction on Neumann boundary $\Gamma_t$, respectively. $\mathbb{E}$ is the constitutive matrix expressing the relationship between force and deformation, in the present work, it is assumed that $\mathbb{E}^i = \ldots = \mathbb{E}^n = \mathbb{E}$ respectively, and it is the same with Poisson's ratio. $\varepsilon$ represents the second order linear strain tensor. $\bar{\boldsymbol{u}}$ is the predefined displacement on Dirichlet boundary $\Gamma_u$, $\boldsymbol{u}$ and $v$ are the displacement field and the potential energy function defined on $\Omega^s = \bigcup_{i=1}^n \Omega_i$ with $u_{ad} = \{ v \mid v \in H^1(\boldsymbol{D}), v = 0 \text{ on } \Gamma_u \}$. $\bar{V}$ is upper bound of the available volume of solid material.

## 3. The framework of ML based parameter tuning strategy

Considering the high performance of modeling of ML, it is taking an increasingly important role in computational mechanics as more data from experiments and calculations. In this study, based on the above introduction, in order to avoid the infeasible optimum and make the MMC based topology optimization approach be free from the manual parameter tuning. An ML based parameter tuning strategy is proposed. It not only saves the manpower but also avoids the local optimum problem caused by the limitations of manual tuning. The highlight of this strategy is that it introduces a data-driven paradigm, which avoids the infeasible optimum solution through the ML based image classifier. The classifier is integrated into the closed loop optimization framework to obtain the reasonably feasible solution. Fig. 5 illustrates the framework of this method.

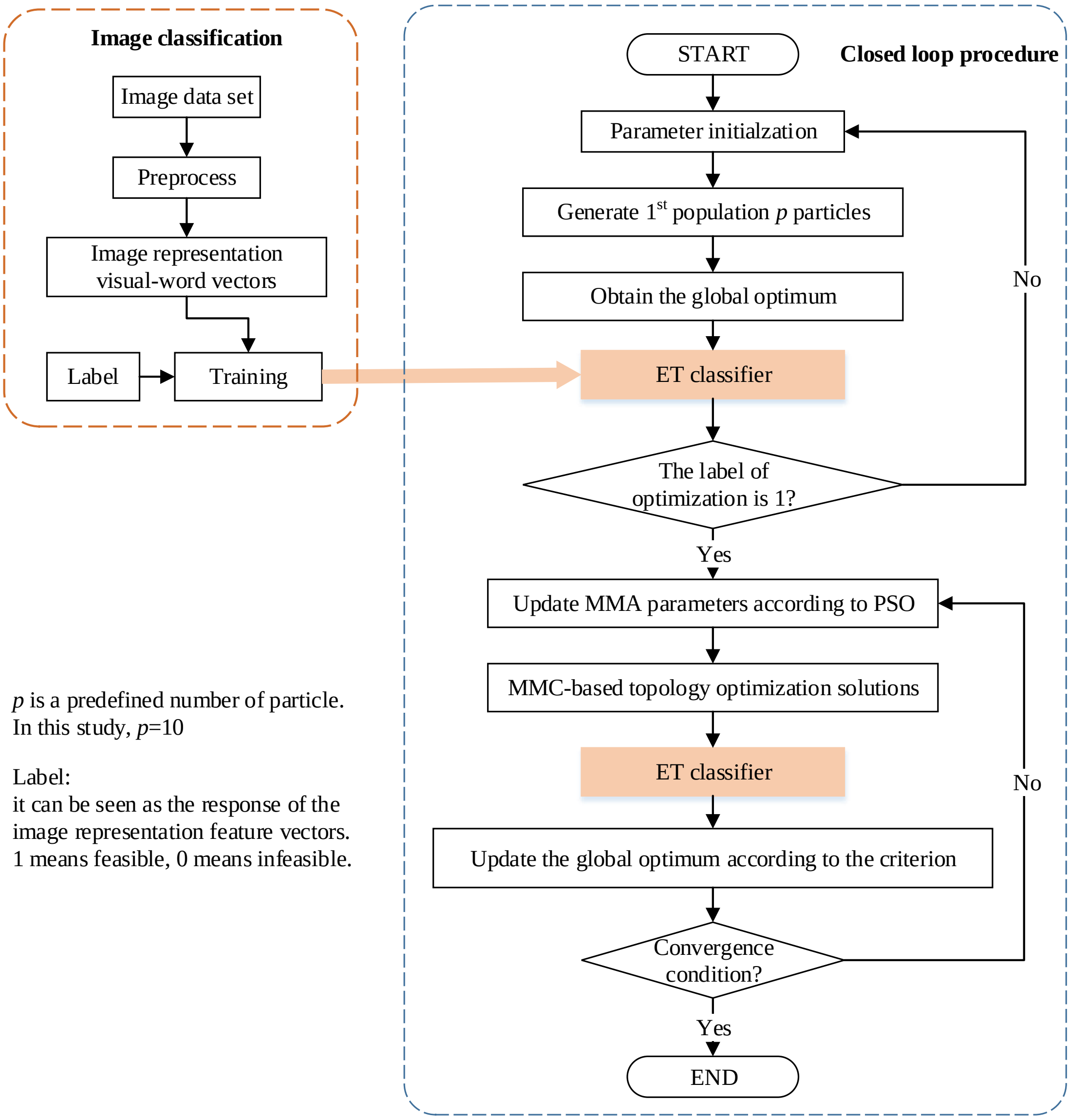

**Fig. 5** Framework of ML based parameter tuning strategy

### 3.1 Image Classification for MMC based topology optimization

In this section, image preprocess and classification will be introduced. Since the initial parameters of MMA of the proposed strategy are randomly selected, the solution with the lowest value of objective function (e.g. *C*) in the initial optimization process might not be feasible. Traditional optimization algorithm (e.g. PSO) might not competent for this task since it might lead the optimization in a wrong direction. Hence, it is particularly important to train a classifier which can effectively determine the feasibility of solutions. Firstly, samples are generated by random MMA parameters in given design space, which are used as the original data for training and testing the classifier subsequently. The data is labeled with the prior knowledge.

### 3.1.1 Image preprocess

It is well-known that describe an image with an effective way in image classification is important. Generally, the selection of a suitable feature space is often based on an actual problem. In this study, since MMC based approach is an explicit algorithm, whose topology optimization solutions are clear in structure and boundary. Therefore, the gradient of pixel value on the geometric boundary is large. These features can easily be detected as key points. These key points are salient image patches that contain rich local information of an image. In the previous work presented in Refs [30, 31], Bag of Visual Words (BoVW) combined with Scale-Invariant Feature Transform (SIFT) [32, 33] has been used as an image representation for image classification well. The details of this procedure are presented as follows.

Firstly, the SIFT is employed to obtain image features through intensive selection of the key points and generates feature descriptors for each key point. This method firstly convolves the image pixel information with the variable-scale Gaussian function to obtain the variable-scale spatial function $L(x,y,\sigma)$ which reflects the local details of the image.

$$L(x,y,\sigma) = G(x,y,\sigma) * I(x,y) \tag{13}$$

$$G(x,y,\sigma) = \frac{1}{2\pi\sigma^2} e^{\frac{-(x^2+y^2)}{2\sigma^2}} \tag{14}$$

where $\sigma$ is the scale space factor, $x, y$ represents the position of key points. $I(x,y)$ Represents the sampling information of the original image under the MMC based approach. More details of the SIFT can be found in Refs [32, 33]. In this way, key points are described by the location, direction, and scale information. As shown in the Fig. 6, all the key points of a topology optimization solution are drawn. According to the gradient value and direction of each key point, the algorithm selects $16\times16$ pixels in the neighborhood and divides them into $4\times4$ small blocks, then statistics the histogram of gradient in 8 directions of each small block. Therefore, descriptors on key points can be represented by $4\times4\times8$ directions, namely a 128-dimensions SIFT feature vectors. Thus, after feature extraction by the SIFT, each image will be described by $\boldsymbol{n}$ 128-dimensional vectors.

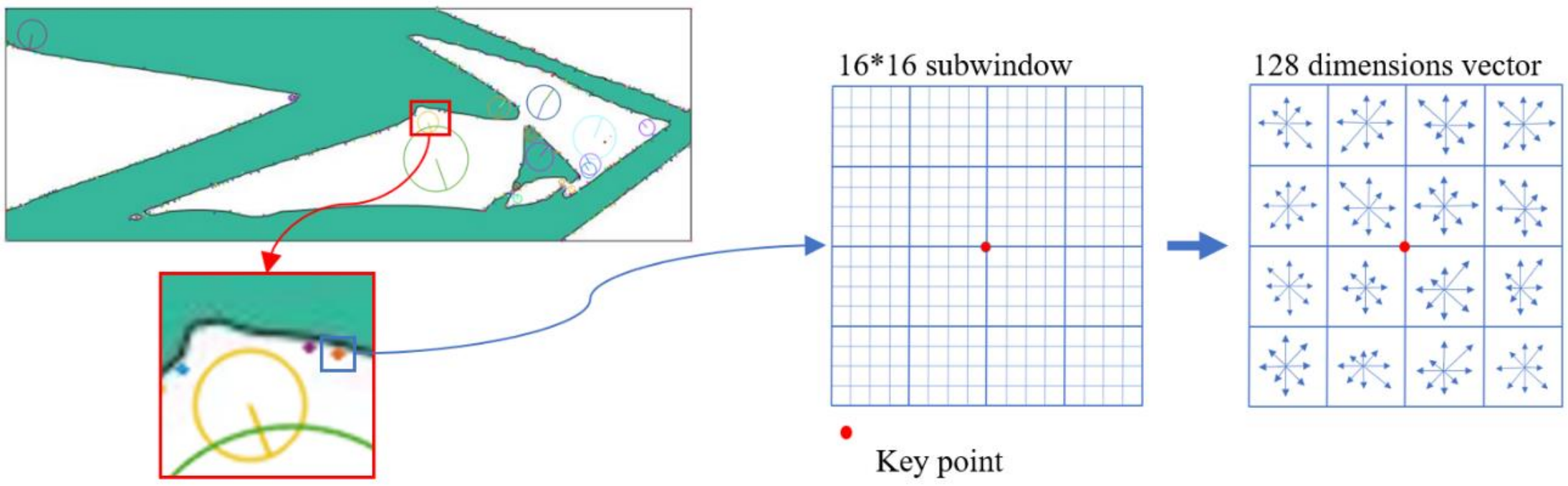

**Fig. 6** The SIFT descriptor of an MMC based topology optimization solution

Secondly, BoVW representation is employed to normalize and quantify all the SIFT feature vectors extracted from the samples. In the last step, images which are generated under the MMC based framework are represented by sets of key point descriptors, but the sets vary in cardinality and lack meaningful ordering. It creates difficulties for ML methods (e.g., classifiers) that require feature vectors of fixed dimension as input. To address this issue, *k*-means [34] clustering algorithm is employed, which clusters the key point descriptors in their feature space into a large number of clusters and encodes each key point by the index of the cluster to which it belongs. Each cluster can be considered as a visual word that represents a specific local pattern shared by the keypoints in that cluster. Hence, the clustering process generates a visual-word vocabulary describing different local patterns in different topology optimization solutions. The number of clusters determines the size of the vocabulary, while there is no consensus as to the appropriate size of a visual word vocabulary, it is acceptable as long as the classification works well. In this study, $k = 64$. Mapping the keypoints to visual words, each topology optimization solution can be described as a "bag of visual words". This representation is analogous to the Bag-of-Words [35] document representation in terms of form and semantics. Both representations are sparse and high dimensional, and just as words convey meanings of a document, visual words reveal local patterns characteristic of the whole image. Finally, as shown in the Fig. 7, converting the "bag of visual words" into a visual-word vector, each MMC topology optimization solution can be described by the frequency of each visual word. In another word it means each image was described and classified by a 64-dimentional visual-word vector which is consist of the frequency of each cluster center.

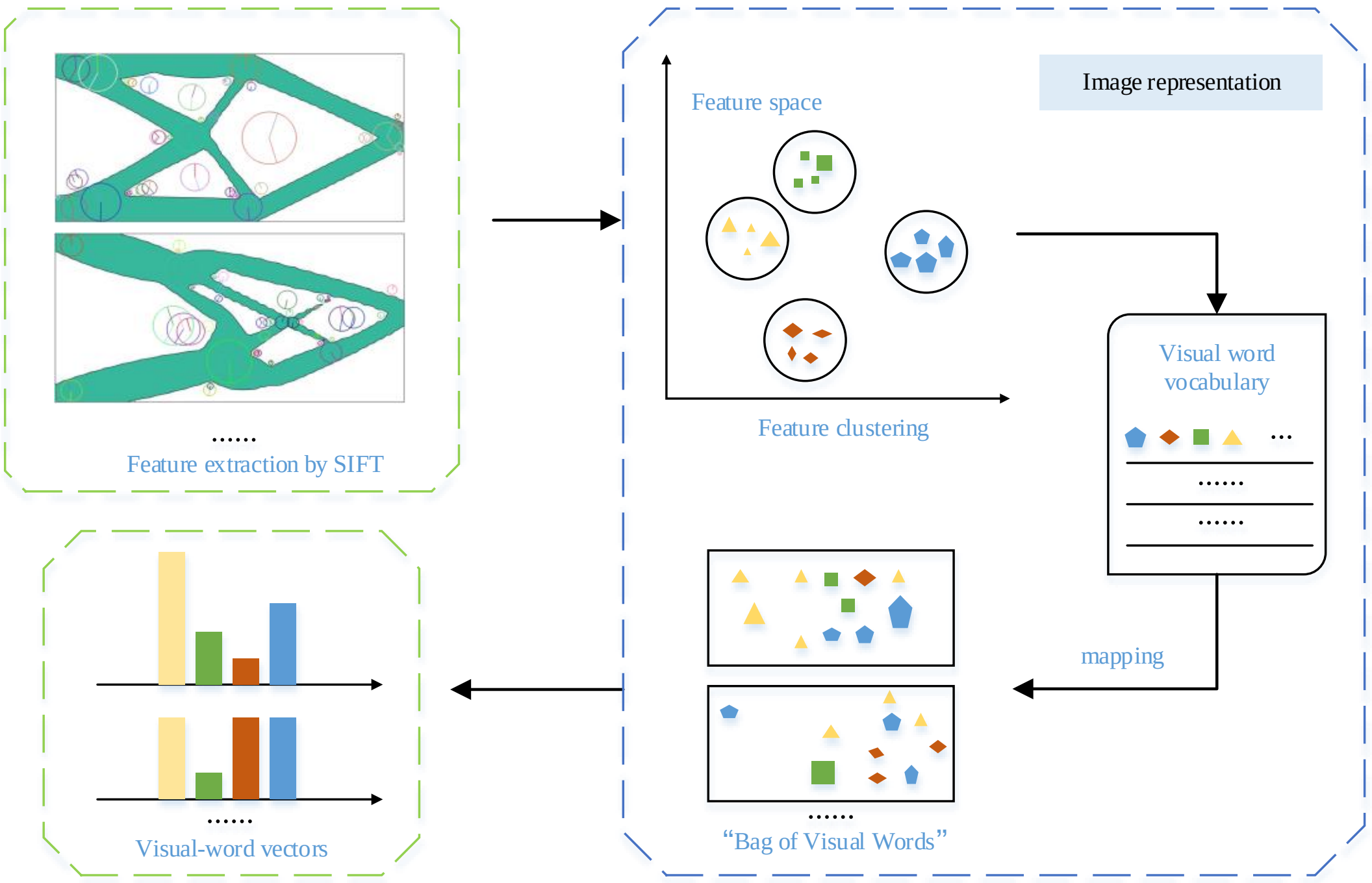

**Fig. 7** Flowchart of image preprocess.

### 3.1.2 Classification based on Extra-Trees (ET)

In the past years, extensive works have been devoted to ML and pattern recognition to construct efficient classifier. Recently, ensemble methods have been of main interest especially in medical image classification [36-38], and these study have proved the efficiency over single classifier on different tasks. Hence, a tree-based ensemble method Extremely Randomized Trees [39] is introduced to MMC based topology optimization solutions' classification in this study. ET builds an ensemble of unpruned decision trees according to the classical top-down procedure. It is an improvement of Random Forest [40]. Compared with other tree-based ensemble methods, The ET splits nodes by choosing cut-points fully at random and uses the whole learning sample (rather than a bootstrap replica) [41] to grow the trees. Therefore, the method randomly picks a single attribute and cut-point at each node, and builds totally randomized trees whose structures are independent of the target variable values of the learning sample. These are also the strengths of this algorithm. From the bias-variance point of view, the rationale behind the ET method is that the explicit randomization of the cut-point and attribute combined with ensemble averaging should be able to reduce variance more strongly. The usage of the full original learning sample rather than bootstrap replicas is motivated in order to minimize bias. These good statistical performances are the reason why the ET algorithm was employed to training the classifier between the visual-word vectors and the labels. The label here can be seen as the classifier's responses to the visual-word vector. In addition, the required amount of training data is much lower than that of deep learning classifier.

Its performance usually depends on three parameters: $K$, the number of attributes randomly selected at each node. $n_{\min}$, the minimum sample size for splitting a node. And $m$, the number of decision trees in an ensemble model. The predictions of the trees are aggregated to yield the final prediction, by majority vote in classification problems and arithmetic average in regression problems. In this study, $K$ and $n_{\min}$ is the default value in order to maximize the computational advantages and autonomy of the method [39]. It is well known that for randomized ensembles the classification error is monotonically decreasing as the increasing of the value of $m$. While, choosing an appropriate value for $m$ is not so important, a high value of m may be selected if the computing ability permits. With many experiments, $m = 400$ can achieve good effects in this study.

### 3.2 Closed loop parameter tuning strategy

Once the classifier is established, it would be integrated into the closed loop optimization framework and the complete parameter tuning strategy is obtained. The classifier provides the feasibility criterion for the parameter tuning, and the optimization is carried out by Particle Swarm Optimization (PSO) algorithm [42-44]. The most common type of implementation defines the particles' behaviors of PSO in two formulas. The first adjusts the velocity or step size of the particle, and the second moves the particle by adding the velocity to its previous position. On each dimension *d*:

$$v_{id}^{t+1} = \alpha_i v_{id}^t + U(0,\beta_1)(p_{id} - x_{id}^t) + U(0,\beta_2)(p_{gd} - x_{id}^t) \tag{15}$$

$$x_{id}^{t+1} = x_{id}^t + v_{id}^{t+1} \tag{16}$$

where *i* is the target particle's index, *d* is the dimension, in this study $d = 4$, $x_i$ is the particle's position that corresponding to the changeable parameters, $v_i$ is the velocity, $p_i$ is the best position found so far by *i*; *g* is the index of *i*'s best neighbor or "global optimum", $\alpha_i$ is a weight factor that decreases with the increase of iteration. $\alpha_{i+1} = 0.99 \times \alpha_i$ and $\beta_1 = 1$, $\beta_2 = 2$ are inertia weight constant and accelerate constant, $U(0,\beta)$ is a uniform random number generator.

When the system works, each particle cycle around a region centered on the centroid of the previous bests $p_i$ and $p_g$. In this method, the optimal parameters update of MMC by PSO not only depend on the value of the objective function, but also satisfy the condition that the corresponding structure of the objective is a feasible structure. Given the random initial parameters in the prescribed design space, the parameters are completely updated by the algorithm until the convergence condition is reached. Specifically, $p_i$ and $p_g$ of PSO will not be an infeasible solution under the proposed approach. The classifier will help the PSO guide the optimization process better. Therefore, in the first particle swarm, it is necessary to obtain a relative feasible solution. Otherwise, the new particle swarm is still generated from the random MMA parameters. The whole framework has been presented in Fig. 5. The PSO is applied due to its simple principle, fewer parameters and easy realization. Moreover, the PSO is memorable, the results of all particles and the corresponding design parameters are saved. It is convenient to analysis the data subsequently.

## 4. Numerical cases

### 4.1 Classification performance criteria

As for the classifier established previously, it is necessary to evaluate the reliability and generalization ability of the model since a reliable model could better guarantee the correct direction of parameter tuning. For binary classification problems, the classification result can be described by a confusion matrix [45] according to the input and output, which is shown in the Table. 1. It can be inferred that the $FP$ means the actual category of the sample is defined as false while the predicted category is positive, and the meanings of $TP$, $TN$ and $FN$ are analogous. Generally, the accuracy of classification is calculated by the Eq. (17)

**Table. 1** Confusion matrix

| - | | **Predicted** | |
|---|---|---|---|
| | | **Positive** | **Negative** |
| **Actual** | **True** | $TP$ | $FN$ |
| | **False** | $FP$ | $TN$ |

$$ACC = \frac{TP + TN}{TP + FN + FP + TN} \tag{17}$$

However, it is difficult for $ACC$ to evaluate classification problems objectively. Hence, Precision ($P$) [46] and recall ($R$) as well as True Positive Rate ($TPR$) and False Positive Rate ($FPR$) are proposed as compensations, which are two pairs of contradictory metrics. $P$ means the proportion of positive actually in all positive samples determined by the classifier, and the $R$ means the proportion of positive determined by the classifier in all actually true samples. $TPR$ and $FPR$ can be inferred from the Eqs. (20) to (21) in the similar way. Moreover, a more comprehensively metrics $F1$ score calculated by Eq. (22) is employed, which is the harmonic mean of the $P$ and $R$. As for these metrics, the larger the better.

$$P = \frac{TP}{TP + FP} \tag{18}$$

$$R = \frac{TP}{TP + FN} \tag{19}$$

$$TPR = \frac{TP}{TP + FN} \tag{20}$$

$$FPR = \frac{FP}{FP + TN} \tag{21}$$

$$F1 = 2 \cdot \frac{P \cdot R}{P + R} \tag{22}$$

In addition to numerical criteria, binary classification problems can also be evaluated according to performance plots. Such as *P-R* plot [47] and Receiver Operating Characteristic (ROC) [48-50] plot. The *P-R* plot intuitively reflects the overall *P*, *R* value of the model, while the ROC plot discusses the generalization ability of the model. The ordinate and abscissa of *P-R* plot are *P* value and *R* value respectively, while the ordinate and abscissa of ROC plot are *TPR* and *FPR*. When compared the properties of different ML methods, if one plot covers the other completely, it can be asserted that the former performs better than the later. While at most cases, there is usually some interlacing in curves. Thus, evaluation criteria need to be considered more comprehensively. Generally, Area Under Curve (*AUC)* is employed to evaluate the model. In this study, by combining *ACC*, $F1$score and *AUC*, the advantages and disadvantages of the model can be comprehensively evaluated.

**4.2 Case I: Cantilever beam**

As shown in Fig. 8, the left side of the short beam is fixed. The length and width of the whole design domain are $DW = 2$ and $DH = 1$, which are discretized into $40 \times 80$ meshes by four-node bilinear elements, and a concentrated vertical load *F* is enforced at the middle point of the right of the beam. The Young's modulus of the material is $\mathbb{E} = 1$ and the Poisson's ratio is $v = 0.3$, the residual volume fraction of the material is 0.4. In this study, there are four changeable parameters whose prescribed design space are given in the Table. 2. Including four adjustable parameters of the MMA. The object of this optimization process is to minimize the Compliance (*C*) of the beam and find the appropriate parameters in the prescribed design space.

**Table. 2** Changeable parameters of Case I

| - | albefa | asyinit | asyincr | asydecr |
|---|---|---|---|---|
| **Upper limit** | 0.75 | 0.75 | 1.50 | 0.75 |
| **Lower limit** | 0.25 | 0.25 | 1.00 | 0.25 |

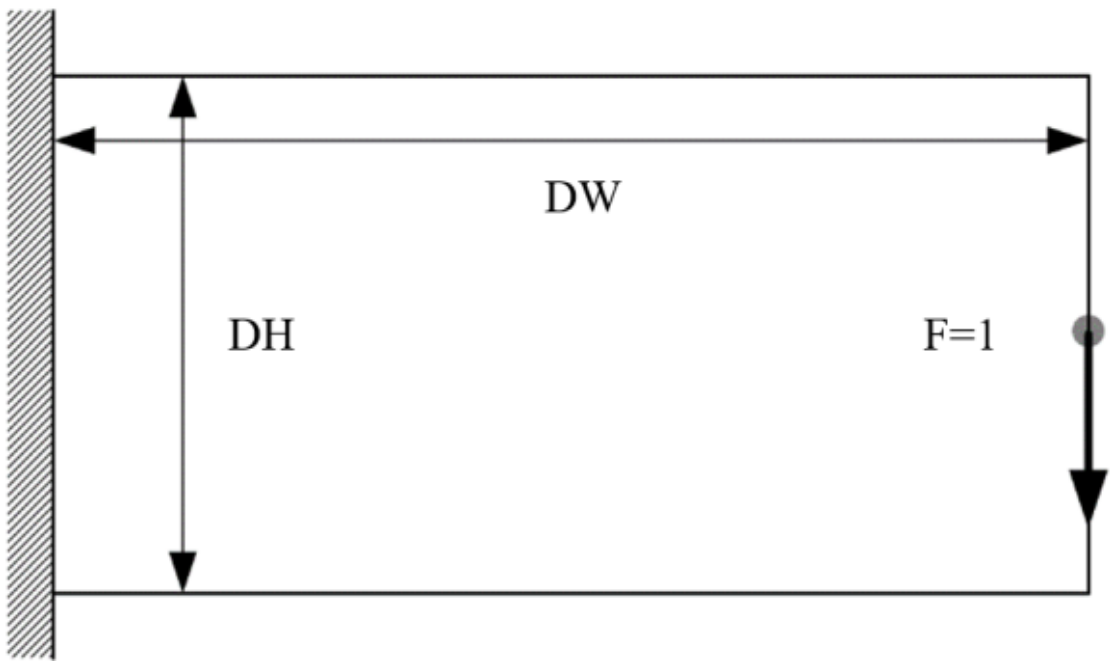

**Fig. 8** Cantilever beam diagram.

In this case, the types of MMC based solutions are shown in the Table. 3. Including the feasible one and the infeasible one. It is obviously that the infeasible structures cannot transfer forces well, and it is not in accordance with the requirements of manufacturing process. As mentioned before, it is caused by the instability of MMC algorithm. Hence, a classifier is trained to avoid the infeasible solution and it assists to guide the process of the parameter tuning. There are total 150 samples in this case, which contains 50 training samples and 100 test samples. In order to highlight the advantages of ET algorithm and guarantee the reliability of the classier, three other classical machine learning algorithms, Random Forest (RF) [51], AdaBoost [52, 53] and GradientBoost [54] are employed to compare in the same sample set.

**Table. 3** Classification standard of Case I

| **Types** | **MMC-based solutions** | | |
|---|---|---|---|
| **Feasible** | 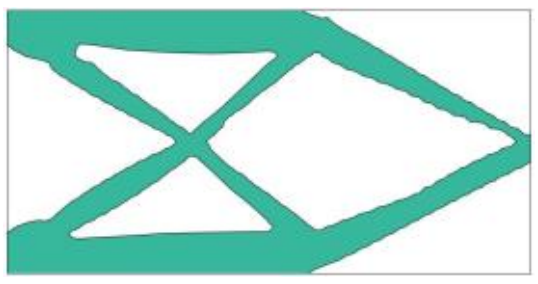 | 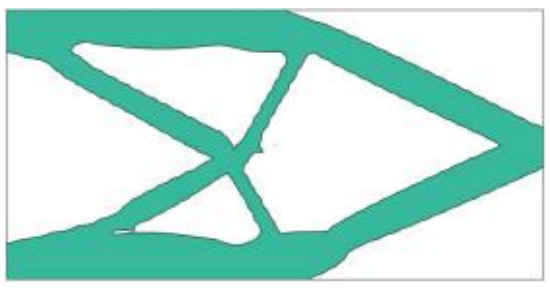 | 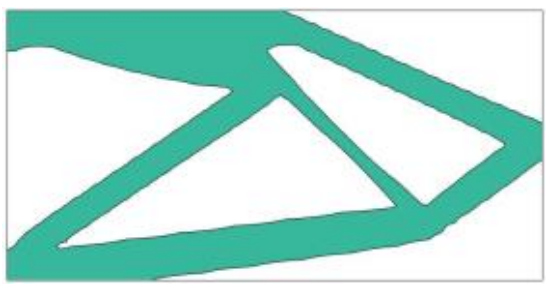 |
| **Infeasible** | 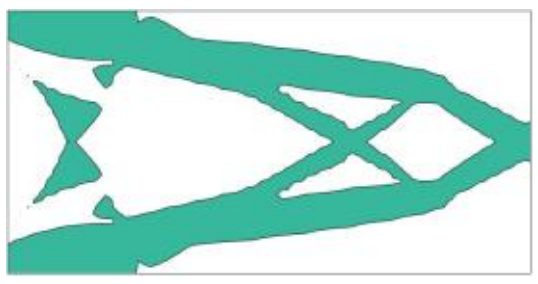 | 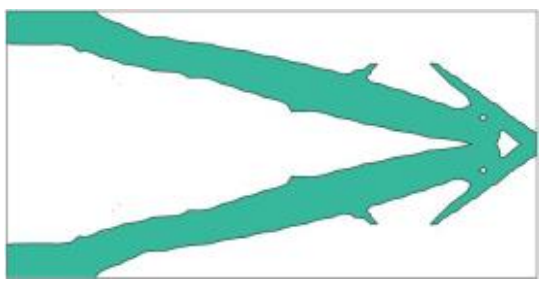 | 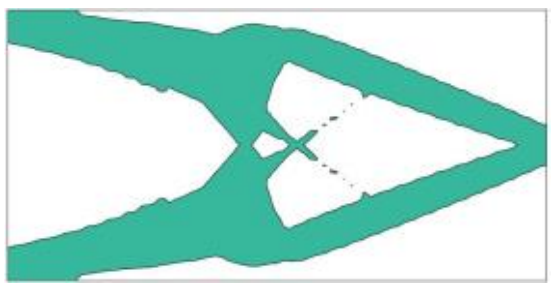 |

The comparison of the *P*-*R* plots and their *AUC* between different ML models is shown in Fig. 9 (a). It can be confirmed that the ET outperforms other algorithms due to its *P*-*R* plot almost completely covers the others. Fig. 9 (b) shows the ROC plot, it also can be confirmed that the most robust model is constructed by the ET, since the *AUC* of the ET (0.96) is the biggest. Moreover, as shown in Table. 4, the ET performs best with the *P* is 0.87, *R* is 0.90, *ACC* is 0.88 and the *F*1 score is 0.88. Thus, it is appropriate to choose ET as the algorithm for the classifier construction. It shows the ET classifier would perform well in the closed loop optimization framework subsequently.

**Table. 4** Performance comparison for 4 algorithms of Case I

| Algorithm | P | R | *F*1 score | ACC |
|---|---|---|---|---|
| **ET** | 0.87 | 0.90 | 0.88 | 0.88 |
| **RF** | 0.81 | 0.88 | 0.85 | 0.84 |
| **AdaBoost** | 0.75 | 0.84 | 0.79 | 0.78 |
| **Gradient Boost** | 0.72 | 0.76 | 0.74 | 0.73 |

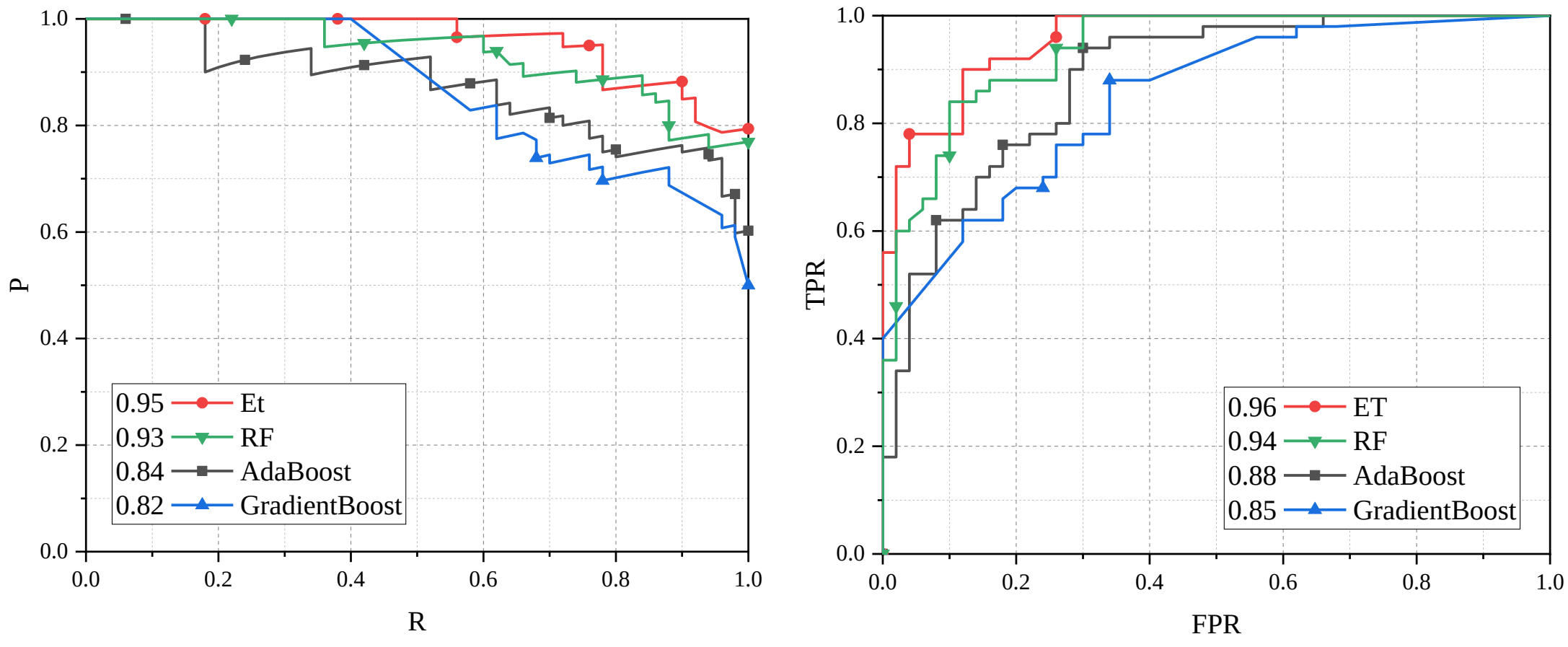

(a) (b)

**Fig. 9** The *P*-*R* plots (a) and the ROC plots (b) of Case I

After the classifier is established, the design parameters for the "global optimal"

structure will be obtained by the PSO algorithm when the optimization process reaches the convergence condition. The convergence curve is shown in the Fig. 10, and the feasible solutions of several pivotal iterations are illustrated as well. It can be found that feasible structures are well obtained from random design parameters by the proposed parameter tuning strategy. The objective function trends to converge after 67$^{th}$ iteration, and the corresponding design parameters of the reasonably feasible solution is shown in Table. 5 with the $C$ is 74.02.

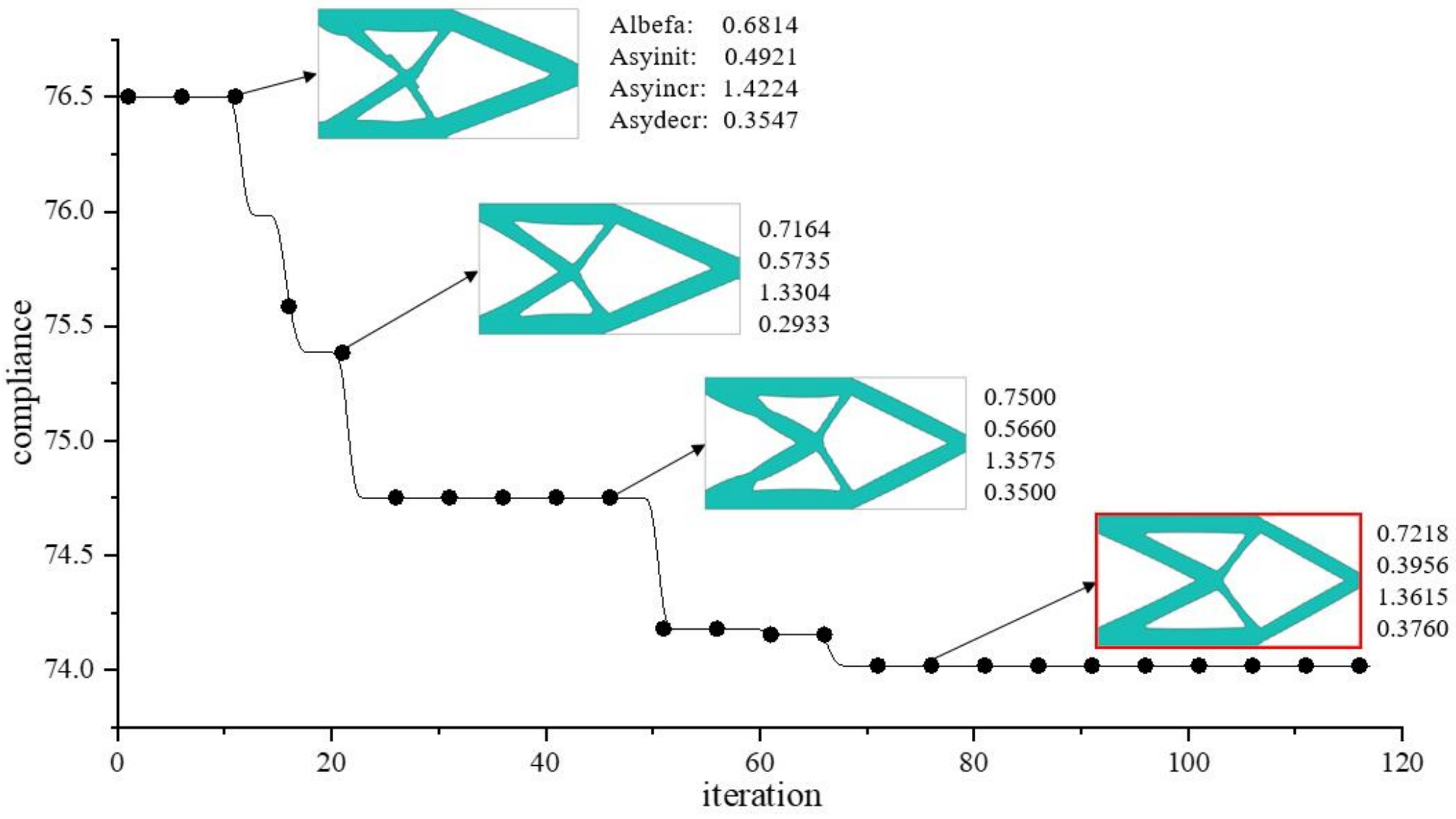

**Fig. 10** The optimization process of the PSO of case I

**Table. 5** Parameters of the reasonably feasible solution of Case I

| **albefa** | **asyinit** | **asyincr** | **asydecr** |
|---|---|---|---|
| 0.7218 | 0.3956 | 1.3615 | 0.3760 |

### 4.2 Case II: L-shape beam

L-shape beam, another mini-compliance case is calculated in this study, whose specific geometry of the design domain and load constraints are depicted in Fig. 11. The top of the beam is fixed and a concentrated vertical downward force $F$ is enforced at the middle of the right of the beam. Moreover, the domain is discretized into 4864 FEM meshes by four-node bilinear elements. As for this case, the objective function, material parameters and the changeable parameters' interval are the same as the Case I.

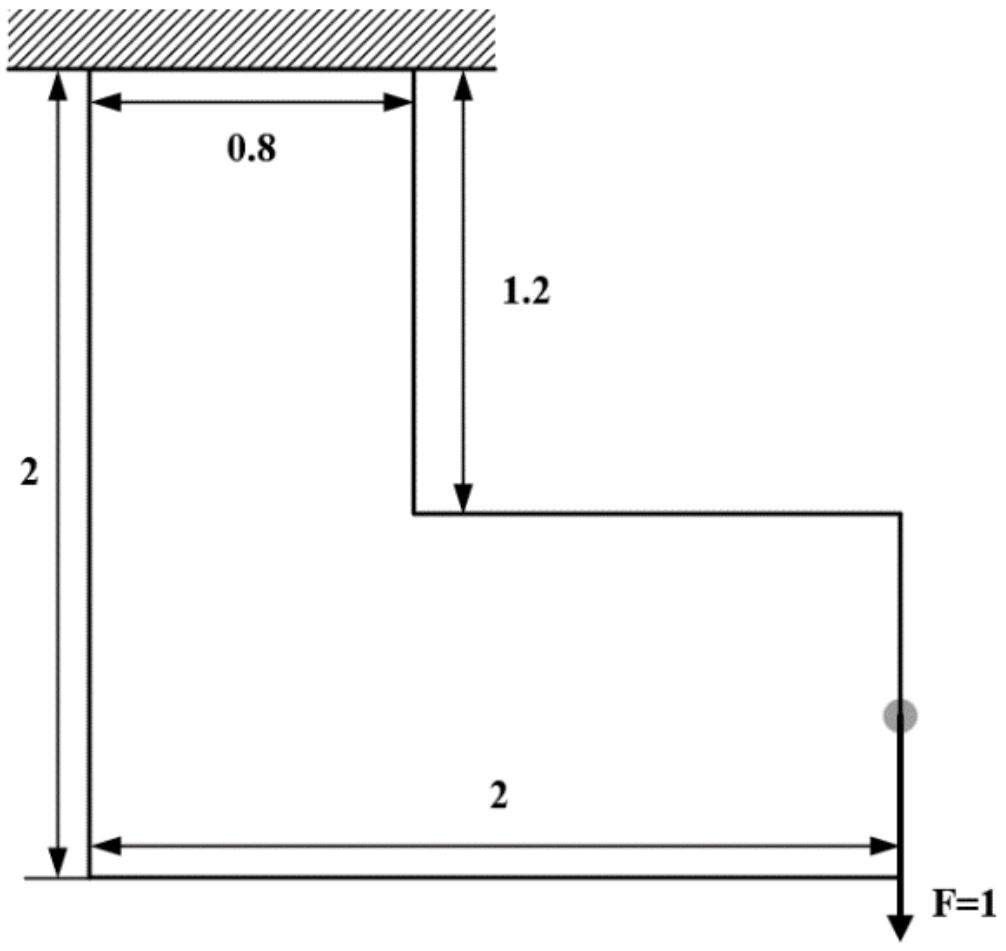

**Fig. 11** L-shape beam diagram

As shown in the Table. 6, similar to the Case I, the solutions of L-shape beam are divided to positive and negative samples. There are also 50 training samples and 100 test samples. The comparison of 4 models is shown in the Table. 7. The *P*-*R* plot and ROC plot are shown in the Fig. 12. Verified by calculation results, the ET algorithm is also believed performs the best with the score *P* is 0.90, *R* is 0.94, *F*1 score is 0.92 and *ACC* is 0.92. In addition, it is also proved that the classification ability of ET is also the most reliable according to the *AUC* of ROC plot and *P*-*R* plot. The ET still maintains its advantages. It is easy for the ET to perform as a good classifier in the optimization process.

**Table. 6** Classification standard of Case II

| Categories | MMC-based solutions | | |
|---|---|---|---|
| **Feasible** | 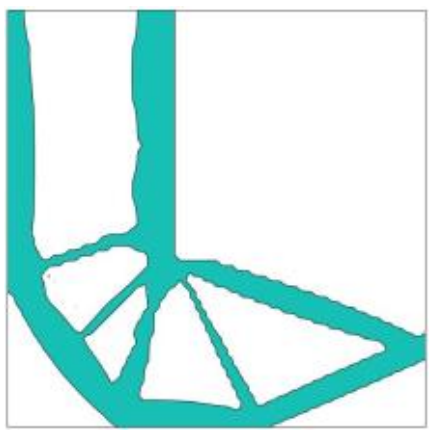 | 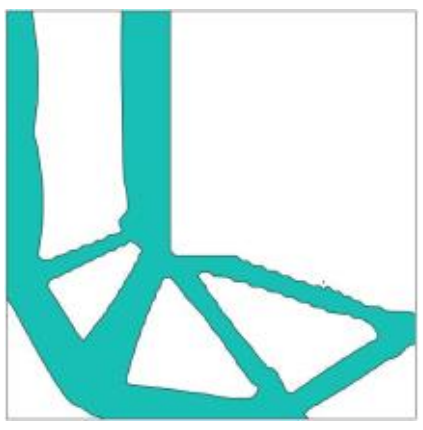 | 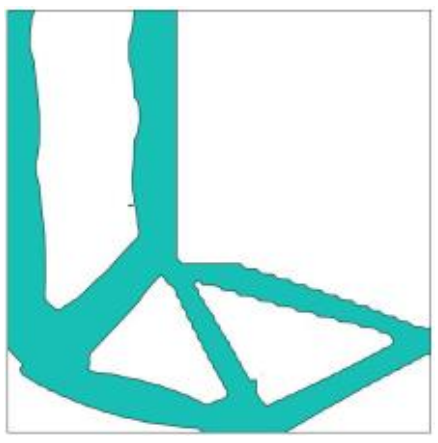 |
| **Infeasible** | 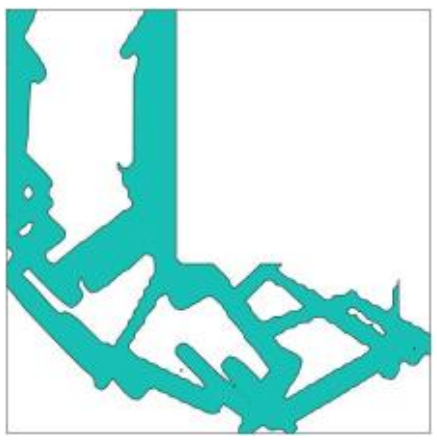 | 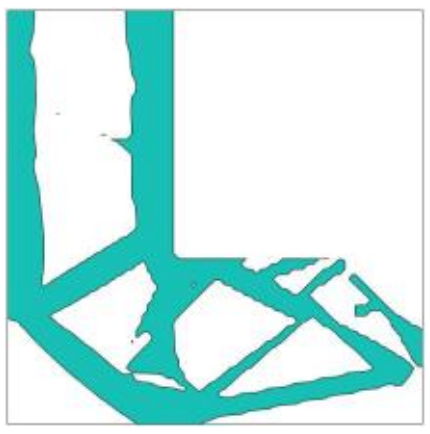 | 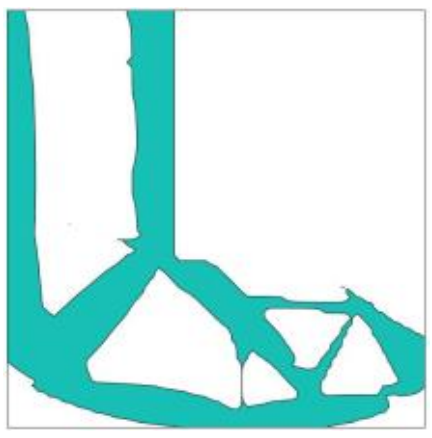 |

**Table. 7** Performance comparison for 4 algorithms of Case II

| Algorithm | P | R | $F1$ score | ACC |
|---|---|---|---|---|
| **ET** | 0.90 | 0.94 | 0.92 | 0.92 |
| **RF** | 0.82 | 0.92 | 0.87 | 0.86 |
| **AdaBoost** | 0.83 | 0.84 | 0.84 | 0.84 |
| **GradientBoost** | 0.83 | 0.78 | 0.80 | 0.81 |

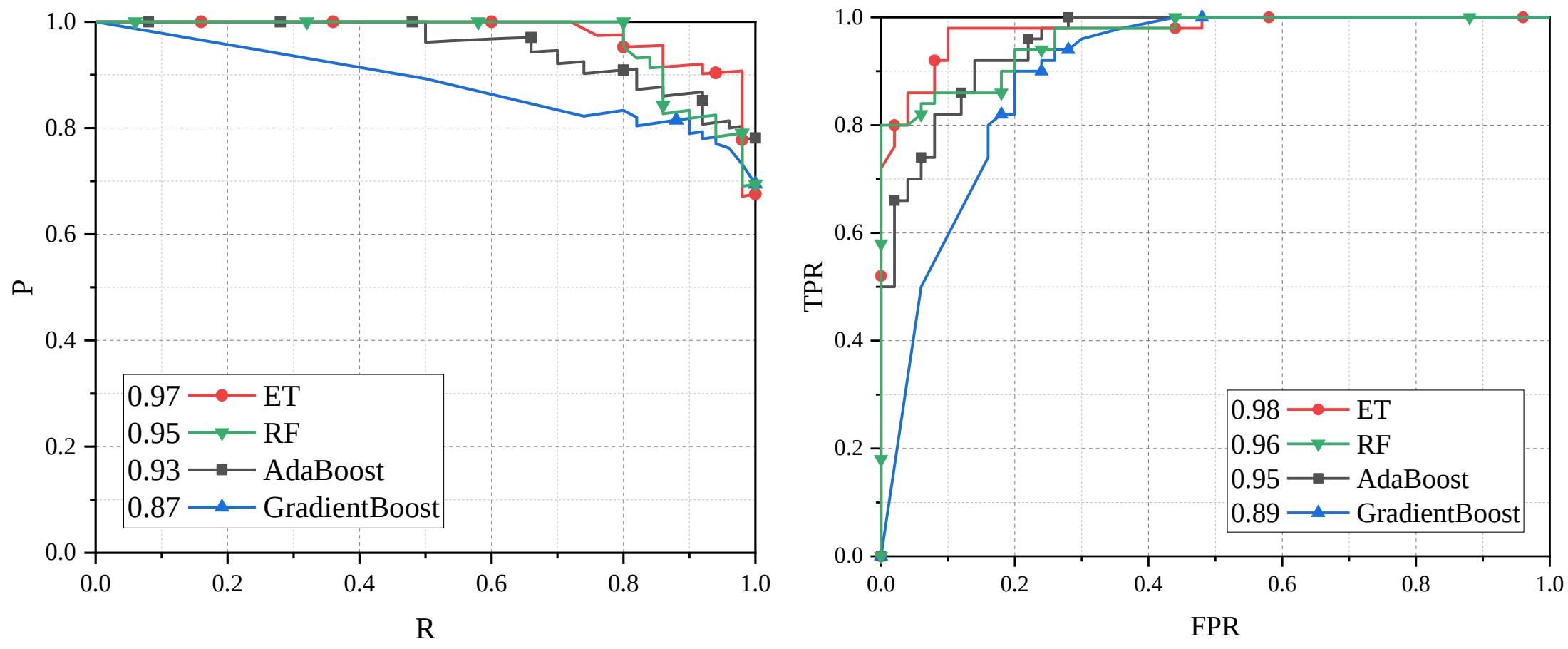

(a) (b)

**Fig. 12** The *P*-*R* plots (a) and the ROC plots (b) of case II

Similar to the Case I, it can be seen from Fig. 13 that the objective function trends to converge as expected after 62th iteration, and the corresponding design parameters of the "global optimal" structure are shown in Table. 8, $C$ is 183.76.

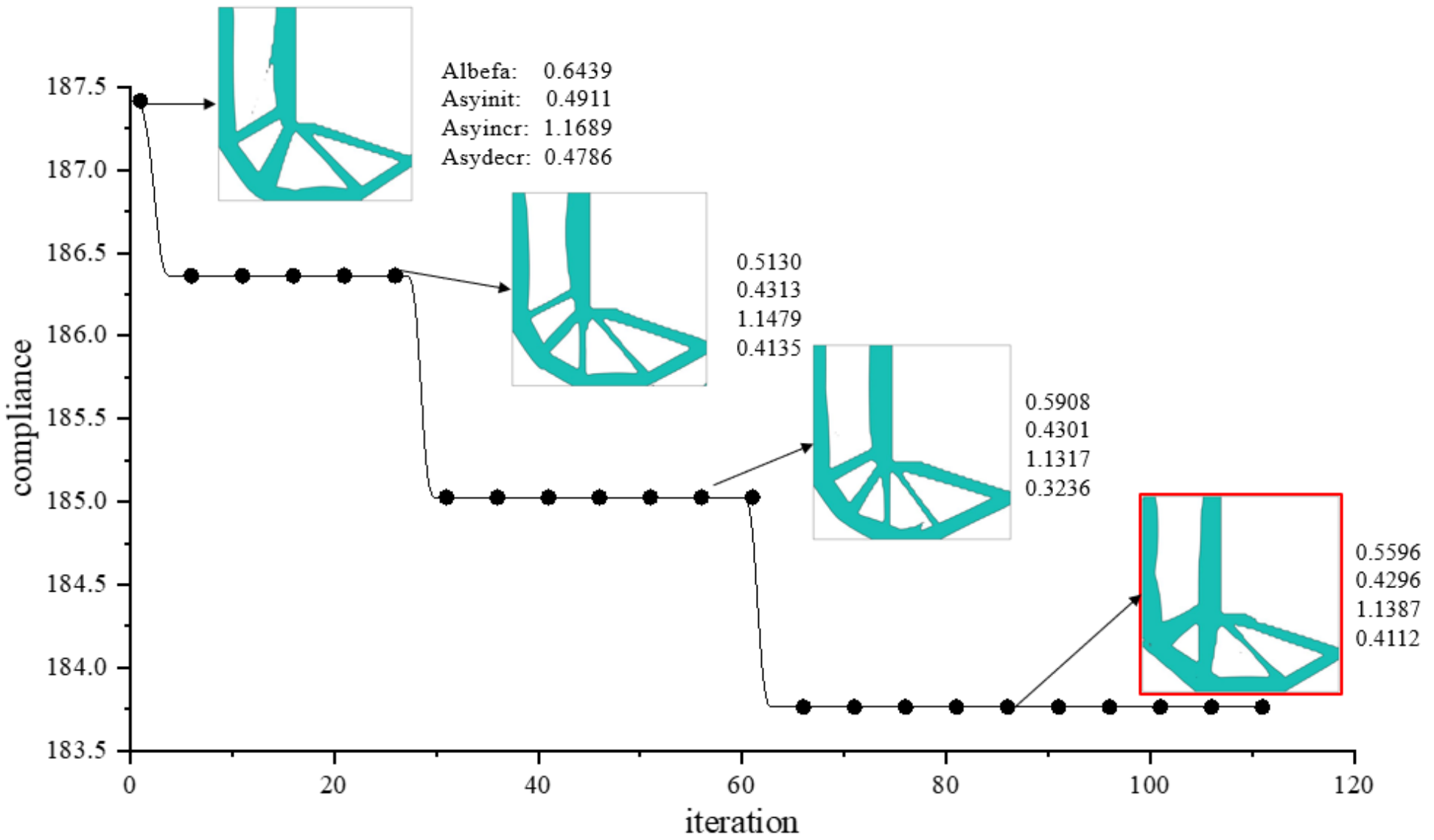

**Fig. 13** The optimization process of the PSO of case II

**Table. 8** Parameters of the reasonably feasible solution of Case II

| albefa | asyinit | asyincr | asydecr |
|---|---|---|---|
| 0.5596 | 0.4296 | 1.1387 | 0.4112 |

## 5. Conclusions

In this study, with the motivation to address the issue of manual tuning of the parameters in MMC based topology optimization. An ML based parameter tuning strategy is proposed. Its highlights can be summarized as follows.

i. As for the MMC based topology optimization, it might obtain many infeasible results due to the inappropriate selection of the parameters in the MMA optimizer. Hence, an ET-based image classifier is introduced to determine whether an optimized solution is feasible or not. It guarantees that the inappropriate parameters would be avoided.

ii. Practically, parameters of MMA are tuned according to the experience and the final obtained optimum might be only an insufficiently feasible solution in the MMC based topology optimization. Two cases are presented to verify the effectiveness of the proposed strategy, it makes the MMC based approach be free from the manual parameter adjustment. It could save a lot of manpower and obtains the reasonably feasible structure in the prescribed design domain.

iii. Moreover, such an ML based parameter tuning strategy can be extended to other fields where the solutions are image-based and the parameters needs a lot of manual tuning.

## Replication of result

All of the code can be found in https://github.com/yoton12138

## Acknowledgments

This work has been supported by Project of the Key Program of National Natural Science Foundation of China under the Grant Numbers 11572120 and 51621004, Key Projects of the Research Foundation of Education Bureau of Hunan Province (17A224).

## Compliance with ethical standards

## Conflict of interest

The authors declare that they have no conflict of interest.

## Reference

1 Sigmund, O., and Maute, K.: 'Topology optimization approaches', Structural and Multidisciplinary Optimization, 2013, 48, (6), pp. 1031-1055

2 Eschenauer, H.A., and Olhoff, N.: 'Topology optimization of continuum structures: A review*', Applied Mechanics Reviews, 2001, 54, (4), pp. 331-390

3 Bendsøe, M.P., and Kikuchi, N.: 'Generating optimal topologies in structural design using a homogenization method', Computer Methods in Applied Mechanics & Engineering, 1988, 71, (2), pp. 197-224

4 Bendsøe, M.P.: 'Optimal shape design as a material distribution problem', Structural optimization, 1989, 1, (4), pp. 193-202

5 Bendsøe, M.P., and Sigmund, O.: 'Material interpolation schemes in topology optimization', Archive of Applied Mechanics, 1999, 69, (9), pp. 635-654

6 Sigmund, O.: 'A 99 line topology optimization code written in Matlab', Structural and Multidisciplinary Optimization, 2001, 21, (2), pp. 120-127

7 Yulin, M., and Xiaoming, W.: 'A level set method for structural topology optimization and its applications', Advances in Engineering Software, 2004, 35, (7), pp. 415-441

8 Wang, M.Y., Wang, X., and Guo, D.: 'A level set method for structural topology optimization', Computer Methods in Applied Mechanics and Engineering, 2003, 192, (1), pp. 227-246

9 Xie, Y.M., and Steven, G.P.: 'A simple evolutionary procedure for structural optimization', Computers & Structures, 1993, 49, (5), pp. 885-896

10 Guo, X., Zhang, W., and Zhong, W.: 'Doing Topology Optimization Explicitly and Geometrically—A New Moving Morphable Components Based Framework', Journal of Applied Mechanics, 2014, 81, (8), pp. 081009-081009-081012

11 Guo, X., Zhang, W., Zhang, J., and Yuan, J.: 'Explicit structural topology optimization based on moving morphable components (MMC) with curved skeletons', Computer Methods in Applied Mechanics and Engineering, 2016, 310, pp. 711-748

12 Zhang, W., Yuan, J., Zhang, J., and Guo, X.: 'A new topology optimization approach based on Moving Morphable Components (MMC) and the ersatz material model', Structural and Multidisciplinary Optimization, 2016, 53, (6), pp. 1243-1260

13 Tcherniak, D.: 'Topology optimization of resonating structures using SIMP method', International Journal for Numerical Methods in Engineering, 2002, 54, (11), pp. 1605-1622

14 Allaire, G., Jouve, F., and Toader, A.-M.: 'Structural optimization using sensitivity analysis and a level-set method', Journal of Computational Physics, 2004, 194, (1), pp. 363-393

15 Sigmund, O.: 'Manufacturing tolerant topology optimization', Acta Mechanica Sinica, 2009, 25, (2), pp. 227-239

16 Guest, J.K., Prévost, J.H., and Belytschko, T.: 'Achieving minimum length scale in topology optimization using nodal design variables and projection functions', International Journal for Numerical

Methods in Engineering, 2004, 61, (2), pp. 238-254
17 Zhang, W., Li, D., Yuan, J., Song, J., and Guo, X.: ‘A new three-dimensional topology optimization method based on moving morphable components (MMCs)’, Computational Mechanics, 2017, 59, (4), pp. 647-665
18 Zhang, W., Chen, J., Zhu, X., Zhou, J., Xue, D., Lei, X., and Guo, X.: ‘Explicit three dimensional topology optimization via Moving Morphable Void (MMV) approach’, Computer Methods in Applied Mechanics and Engineering, 2017, 322, pp. 590-614
19 Zhou, M., Lazarov, B.S., Wang, F., and Sigmund, O.: ‘Minimum length scale in topology optimization by geometric constraints’, Computer Methods in Applied Mechanics and Engineering, 2015, 293, pp. 266-282
20 Zhang, W., Li, D., Zhou, J., Du, Z., Li, B., and Guo, X.: ‘A Moving Morphable Void (MMV)-based explicit approach for topology optimization considering stress constraints’, Computer Methods in Applied Mechanics and Engineering, 2018, 334, pp. 381-413
21 Svanberg, K.: ‘The method of moving asymptotes—a new method for structural optimization’, International Journal for Numerical Methods in Engineering, 1987, 24, (2), pp. 359-373
22 Shen, S.R.M.-H.H.: ‘Neural networks for topology optimization’, Machine Learning, 2018.8, (Department of Mechanical and Aerospace Engineering)
23 Sosnovik I, O.I.: ‘A novel topology design approach using an integrated deep learning network architecture’, 2017, (Department of Mechanical and Aerospace Engineering)
24 Banga, S.: ‘3D Topology Optimization Using Convolutional Neural Networks’, Machine Learning, 2018, (Carnegie Mellon University )
25 Lei, X., Liu, C., Du, Z., Zhang, W., and Guo, X.: ‘Machine Learning-Driven Real-Time Topology Optimization Under Moving Morphable Component-Based Framework’, Journal of Applied Mechanics, 2018, 86, (1), pp. 011004-011004-011009
26 Amstutz, S., and Andrä, H.: ‘A new algorithm for topology optimization using a level-set method’, Journal of Computational Physics, 2006, 216, (2), pp. 573-588
27 Hoang, V.-N., and Jang, G.-W.: ‘Topology optimization using moving morphable bars for versatile thickness control’, Computer Methods in Applied Mechanics and Engineering, 2017, 317, pp. 153-173
28 Chen, S., Wang, M.Y., and Liu, A.Q.: ‘Shape feature control in structural topology optimization’, Computer-Aided Design, 2008, 40, (9), pp. 951-962
29 Guest, J.K., Asadpoure, A., and Ha, S.-H.: ‘Eliminating beta-continuation from Heaviside projection and density filter algorithms’, Structural and Multidisciplinary Optimization, 2011, 44, (4), pp. 443-453
30 Shekhar, R., and Jawahar, C.V.: ‘Word Image Retrieval Using Bag of Visual Words’, in Editor (Ed.)^(Eds.): ‘Book Word Image Retrieval Using Bag of Visual Words’ (2012, edn.), pp. 297-301
31 Yang, J., Jiang, Y.-G., Hauptmann, A.G., and Ngo, C.-W.: ‘Evaluating bag-of-visual-words representations in scene classification’. Proc. Proceedings of the international workshop on Workshop on multimedia information retrieval, Augsburg, Bavaria, Germany2007 pp. Pages
32 Lowe, D.G.: ‘Object Recognition from Local Scale-Invariant Features’, 1999, (University of British Columbia)
33 Cruz-Mota, J., Bogdanova, I., Paquier, B., Bierlaire, M., and Thiran, J.-P.: ‘Scale Invariant Feature Transform on the Sphere: Theory and Applications’, International Journal of Computer Vision, 2012, 98, (2), pp. 217-241
34 Hartigan, J.A., and Wong, M.A.: ‘Algorithm AS 136: A K-Means Clustering Algorithm’, Journal of the Royal Statistical Society. Series C (Applied Statistics), 1979, 28, (1), pp. 100-108
35 Wallach, H.M.: ‘Topic modeling: beyond bag-of-words’. Proc. Proceedings of the 23rd international conference on Machine learning, Pittsburgh, Pennsylvania, USA2006 pp. Pages
36 Desir, C., Petitjean, C., Heutte, L., Salaun, M., and Thiberville, L.: ‘Classification of Endomicroscopic Images of the Lung Based on Random Subwindows and Extra-Trees’, IEEE Transactions on Biomedical Engineering, 2012, 59, (9), pp. 2677-2683
37 Marée, R., Geurts, P., and Wehenkel, L.: ‘Random subwindows and extremely randomized trees for image classification in cell biology’, Bmc Cell Biology, 2007, 8, (1), pp. 1-12
38 Marée, R., Geurts, P., and Wehenkel, L.: ‘Biological Image Classification with Random Subwindows and Extra-Trees’, English, 2006, 3, (1), pp. 75-90
39 Geurts, P., Ernst, D., and Wehenkel, L.: ‘Extremely randomized trees’, Machine Learning, 2006, 63, (1), pp. 3-42
40 Breiman, L.: ‘Random Forests’, Machine Learning, 2001, 45, (1), pp. 5-32
41 Lemmens, A., and Croux, C.: ‘Bagging and Boosting Classification Trees to Predict Churn’, Journal of Marketing Research, 2006, 43, (2), pp. 276-286
42 Eberhart, and Yuhui, S.: ‘Particle swarm optimization: developments, applications and resources’, in Editor (Ed.)^(Eds.): ‘Book Particle swarm optimization: developments, applications and resources’ (2001, edn.), pp. 81-86 vol. 81
43 Kennedy, J., and Eberhart, R.: ‘Particle swarm optimization’, in Editor (Ed.)^(Eds.): ‘Book Particle swarm optimization’ (2002, edn.), pp.
44 Kannan, G., Noorul Haq, A., and Devika, M.: ‘Analysis of closed loop supply chain using genetic algorithm and particle swarm optimisation’, International Journal of Production Research, 2009, 47, (5), pp. 1175-1200
45 Townsend, J.T.: ‘Theoretical analysis of an alphabetic confusion matrix’, Perception & Psychophysics, 1971, 9, (1), pp. 40-50
46 Goutte, C., and Gaussier, E.: ‘A Probabilistic Interpretation of Precision, Recall and F-Score, with Implication for Evaluation’, in Editor (Ed.)^(Eds.): ‘Book A Probabilistic Interpretation of Precision, Recall and F-Score, with Implication for Evaluation’ (Springer Berlin Heidelberg, 2005, edn.), pp. 345-359
47 Saito, T., and Rehmsmeier, M.: ‘The Precision-Recall Plot Is More Informative than the ROC Plot When Evaluating Binary Classifiers on Imbalanced Datasets’, PLOS ONE, 2015, 10, (3), pp. e0118432
48 Hanley, J.A., and McNeil, B.J.: ‘The meaning and use of the area under a receiver operating characteristic (ROC) curve’, Radiology, 1982, 143, (1), pp. 29-36
49 Hand, D.J., and Till, R.J.: ‘A Simple Generalisation of the Area Under the ROC Curve for Multiple Class Classification Problems’, Machine Learning, 2001, 45, (2), pp. 171-186
50 Davis, J., and Goadrich, M.: ‘The relationship between Precision-Recall and ROC curves’. Proc. Proceedings of the 23rd international conference on Machine learning, Pittsburgh, Pennsylvania, USA2006 pp. Pages
51 Liaw, A., and Wiener, M.: ‘Classification and Regression by RandomForest’ (2001. 2001)
52 Collins, M., Schapire, R.E., and Singer, Y.: ‘Logistic Regression, AdaBoost and Bregman Distances’, Machine Learning, 2002, 48, (1), pp. 253-285
53 Khammari, A., Nashashibi, F., Abramson, Y., and Laurgeau, C.: ‘Vehicle detection combining gradient analysis and AdaBoost classification’, in Editor (Ed.)^(Eds.): ‘Book Vehicle detection combining gradient analysis and AdaBoost classification’ (2005, edn.), pp. 66-71 54 Lawrence, R., Bunn, A., Powell, S., and Zambon, M.: ‘Classification of remotely sensed imagery using stochastic gradient boosting as a refinement of classification tree analysis’, Remote Sensing of Environment, 2004, 90, (3), pp. 331-336